\documentclass[pdflatex,sn-mathphys-num]{sn-jnl}

\usepackage{graphicx}
\usepackage{amsmath,amssymb,amsfonts}
\usepackage{amsthm}
\usepackage{mathrsfs}
\usepackage{enumitem}
\usepackage[title]{appendix}
\usepackage{xcolor}
\usepackage{manyfoot}
\usepackage{booktabs}
\usepackage{listings}

\theoremstyle{thmstyleone}
\newtheorem{theorem}{Theorem}
\newtheorem{lemma}[theorem]{Lemma}
\newtheorem{corollary}[theorem]{Corollary}
\newtheorem{proposition}[theorem]{Proposition} 

\theoremstyle{thmstyletwo}

\newtheorem{remark}{Remark}

\theoremstyle{thmstylethree}

\newtheorem{hypothesis}{Hypothesis}[section]

\DeclareMathOperator*{\argmin}{arg\,min}

\DeclareMathOperator{\prox}{prox}
\newcommand{\qbox}[1]{\quad\hbox{#1}\quad}

\newcommand{\R}{\mathbb{R}}

\newcommand{\ackname}{Acknowledgements}
\def\acknowledgement{\par\addvspace{17pt}\small\rmfamily
\trivlist\if!\ackname!\item[]\else
\item[\hskip\labelsep
{\bfseries\ackname}]\fi}

\newenvironment{acknowledgements}{\begin{acknowledgement}}
{\end{acknowledgement}}

\begin{document}

\title{ Convergence Rate Analysis for Monotone Accelerated Proximal Gradient Method }

\author{\fnm{Zepeng} \sur{Wang}}\email{zepeng.wang@rug.nl}

\author{\fnm{Juan} \sur{Peypouquet}}\email{j.g.peypouquet@rug.nl}

\affil{\orgdiv{Bernoulli Institute for Mathematics, Computer Science and Artificial Intelligence}, \orgname{University of Groningen}, \state{Groningen}, \country{The Netherlands}}

\abstract{We propose a monotone accelerated proximal gradient method for solving convex composite optimization problems, guaranteeing nonincreasing function values along the iterates — a property that improves numerical stability and is not enjoyed by standard accelerated schemes. The method generalizes the Monotone FISTA algorithm of Beck and Teboulle. In the convex setting, we establish the optimal rate of $\mathcal{O}\left( \frac{1}{k^2} \right)$ and show that all weak subsequential limit points of the iterates are minimizers. In the strongly convex setting, we establish a linear rate $\mathcal{O}\left( \frac{1}{k^2}(1+\rho)^{-k} \right)$ with $\rho\sim \frac{\mu}{3L}$, without requiring knowledge of the strong convexity parameter — roughly a five-fold improvement in the contraction constant over the best previously known rate for monotone methods, and more than 30\% larger than that of the best known rate for non-monotone accelerated methods. Following a similar idea, we propose a variant of Nesterov's accelerated proximal method and establish a linear rate under strong convexity, which is the same as the one above, and is faster than the known results for non-monotone methods.}



\maketitle

\section{Introduction}
Let $H$ be a real Hilbert space, and consider the convex composite optimization problem:
\begin{equation}\label{Prob: cvx_composite}
\min_{x\in H} F(x) = f(x) + g(x),
\end{equation}
where $f: H\to\mathbb{R}$ is convex and $L$-smooth, and $g:H\to\mathbb{R}\cup\{\infty\}$ is convex, proper and lower-semicontinuous. We assume that the solution set of \eqref{Prob: cvx_composite} is nonempty, and write $x^*\in\argmin(F)$ and $F^* = F(x^*)$.

To solve \eqref{Prob: cvx_composite}, one can employ Nesterov's accelerated proximal method, also referred to as FISTA \cite{Beck_2009}, which was developed based on Nesterov's acceleration technique \cite{Nesterov_1983}. It takes the form \cite{Beck_2009}: 
\begin{equation}\label{Algo: APM}\tag{APM}
\left\{
\begin{aligned}
y_{k+1} &= \prox_{sg}\left( x_k - s\nabla f(x_k) \right),\\
x_{k+1} &= y_{k+1} + \frac{t_k-1}{t_{k+1}}(y_{k+1}-y_k),
\end{aligned}
\right.
\end{equation}
where $s\in \left( 0, 1/L \right]$ and
\begin{equation}\label{E: t_k}
t_k = \left\{
\begin{array}{ccl}
1,&\text{if}& k=0,\\[3pt]
\frac{1+\sqrt{1+4t_{k-1}^2}}{2},&\text{if}& k\ge 1.
\end{array}
\right.
\end{equation} 
In practice, $t_k = \frac{k+\alpha-1}{\alpha-1}$ is a common choice such that the extrapolation parameter $\frac{t_k-1}{t_{k+1}} = \frac{k}{k+\alpha}$ with $\alpha>1$. For $\alpha\ge 3$, a convergence rate of $\mathcal{O}\left( \frac{1}{k^2} \right)$ for the function values was shown in \cite{Beck_2009}. But one can actually obtain a faster rate $o\left( \frac{1}{k^2} \right)$ with $\alpha > 3$ \cite{Attouch_2016}. When $f$ is convex and satisfies a local error-bound condition, a convergence rate of $o\left( \frac{1}{k^{2\alpha}} \right)$ is guaranteed, as long as $\alpha>1$ and $s\in(0,1/L)$ \cite{Liu_2024}. When $f$ is $\mu$-strongly convex, a linear convergence rate $\mathcal{O}\left( \frac{1}{k^2}(1 + \rho)^{-k}\right)$, with $\rho=\frac{\mu}{16L}$, holds for $\alpha\ge 3$ and $s=\frac{1}{2L}$ \cite{Shi_2024,Shi_2024_Monotone}. This rate was further improved to $\mathcal{O}\left( \frac{1}{k^2}(1+\rho)^{-k}\right)$, with $\rho=\frac{\mu}{4L}$ \cite{Bao_2023}. This linear convergence rate is achieved without knowing the strong convexity parameter $\mu$, which may be difficult or computationally expensive to estimate in practice. If $\mu$ is known, an accelerated linear convergence rate $\mathcal{O}\left( \left(1 - \rho \right)^k \right)$, with $\rho=\sqrt{\mu/L}$ is obtained by replacing the extrapolation parameter $\frac{t_k-1}{t_{k+1}}$ in \eqref{Algo: APM} by $\frac{1-\sqrt{\mu/L}}{1+\sqrt{\mu/L}}$ \cite{Nesterov_2004}.

Despite their fast convergence, the function values on the sequences generated by \eqref{Algo: APM} are, in general, {\it not nonincreasing}. In practice, this complicates monitoring convergence and may undermine performance. The lack of monotonicity may also result in divergence when \eqref{Algo: APM} involves inexact computations  \cite[Section V.A]{Beck_2009_Monotone}. This problem can be fixed by using restarting techniques \cite{Boyd_2014,Candes_2015,Alamo_2023,Juan_2023,Aujol_2024}, for example, but we shall not pursue this line of research here. Another solution is to force the monotonicity of the function values by structurally modifying \eqref{Algo: APM}. This is achieved by the monotone accelerated proximal method, based on \cite{Beck_2009_Monotone}:
\begin{equation}\label{Algo: M-APM}\tag{M-APM}
\left\{
\begin{aligned}
z_k &= \prox_{sg}\left( x_k - s\nabla f(x_k) \right),\\
y_{k+1} &= \left\{
\begin{array}{ccl}
z_k,&&\text{if  }F(z_k)\le F(y_k),\\[3pt]
y_k,&&\text{otherwise},
\end{array}
\right.\\
x_{k+1} &= y_{k+1} + \frac{t_k-1}{t_{k+1}}(y_{k+1}-y_k) 
           + \frac{t_k}{t_{k+1}}(z_k-y_{k+1})
           + \frac{t_k(\gamma-1)}{t_{k+1}}(z_k-x_k),
\end{aligned}
\right.
\end{equation}
where $s\in \left( 0, 1/L \right]$, $\gamma\in(0,\infty)$ and $(t_k)$ is given by \eqref{E: t_k}. When $\gamma=1$, \eqref{Algo: M-APM} reduces to the monotone FISTA (M-FISTA) algorithm proposed in \cite{Beck_2009_Monotone}. The convergence rate of M-FISTA was studied in \cite{Beck_2009_Monotone} and recently in \cite{Shi_2024_Monotone}. In particular, for M-FISTA, one obtains a sublinear rate $\mathcal{O}\left( \frac{1}{k^2} \right)$ in the convex setting \cite{Beck_2009_Monotone}, and a linear rate $\mathcal{O}\left( \frac{1}{k^2}(1 + \rho)^{-k}\right)$ with $\rho=\frac{\mu}{16L}$ for $s=1/(2L)$ in case $f$ is strongly convex \cite{Shi_2024_Monotone}. 

Compared to M-FISTA \cite{Beck_2009_Monotone}, \eqref{Algo: M-APM} features the inclusion of a correction term:
\begin{equation}\label{E: correction_term}
t_k/t_{k+1}(\gamma-1)(z_k-x_k).
\end{equation}
The motivation for including this term comes from a factor-$\sqrt{2}$ speedup of Optimized Gradient Method (OGM) \cite{Kim_2016,Kim_2017,Park_2023} upon Nesterov's accelerated gradient method (NAG) \cite{Nesterov_1983}. Note that OGM and NAG can be written in a unifying form: 
\begin{equation*}
\left\{
\begin{aligned}
y_{k+1} &= x_k - s\nabla f(x_k),\\
x_{k+1} &= y_{k+1} + \frac{t_k-1}{t_{k+1}}(y_{k+1}-y_k) + \frac{t_k(\gamma-1)}{t_{k+1}}(y_{k+1}-x_k),
\end{aligned}
\right.
\end{equation*}
which reduces to OGM when $\gamma=2$ and NAG when $\gamma=1$. We introduce the parameter $\gamma$ in \eqref{Algo: M-APM} so that we have more control on manipulating the convergence behavior of \eqref{Algo: M-APM}. This turns out to be helpful in deriving a faster convergence rate for \eqref{Algo: M-APM} under strong convexity. The correction term motivates the energy sequence used in the convergence proof (see Section \ref{Sec: energy_estimations}). Tuning $\gamma$ lets us trade off the correction term against the extrapolation term, which in the end yields the improved contraction constant $\rho\sim\mu/3L$ in Section \ref{Sec: strongly convex}. 

Our contributions are summarized as follows and comparisons with existing results are presented in Table \ref{Tab: comparison}.

\begin{itemize}
\item 
We establish a fast convergence rate $\mathcal{O}\left( \frac{1}{k^2} \right)$ for \eqref{Algo: M-APM} in the convex setting. As a byproduct, all weak subsequential limit points of the iterates are minimizers of the objective function.

\item
In the strongly convex setting, we derive a linear convergence rate of $\mathcal{O}\left( \frac{1}{k^2}(1+\rho)^{-k} \right)$, with $\rho\sim \frac{\mu}{3L}$, for \eqref{Algo: M-APM}. This means roughly a five-fold increase in the contraction constant, and hence a substantially faster linear rate, improving upon the result in \cite{Shi_2024_Monotone}, with $\rho\sim\frac{\mu}{16L}$. The improvement is attributed to the correction term \eqref{E: correction_term}.

\item Interestingly, the contraction constant found above is also more than 30\% larger than that of the best known rate for \eqref{Algo: APM} \cite{Bao_2023}. Following the same idea, we then propose a variant of \eqref{Algo: APM} by the inclusion of a similar correction term, and establish a linear rate in the strongly convex setting, which is the same as that for \eqref{Algo: M-APM}, thus is faster than the known results for \eqref{Algo: APM} and its variants \cite{Shi_2024,Shi_2024_Monotone,Bao_2023}.
\end{itemize}

\begin{table}[ht]
\centering
\caption{Linear-rate constants $\rho$ for accelerated methods under strong convexity, requiring no knowledge of $\mu$ (excluding \cite{Nesterov_2004}, which does require $\mu$).}
\label{Tab: comparison}
\begin{tabular}{@{}llll@{}}
\toprule
Method & Step size $s$ & Rate constant $\rho$ & Reference \\
\midrule
Non-monotone \eqref{Algo: APM} & $1/2L$ & $\mu/16L$ & \cite{Shi_2024,Shi_2024_Monotone} \\
Non-monotone \eqref{Algo: APM} & $1/2L$ & $\mu/4L$ & \cite{Bao_2023} \\
Non-monotone \eqref{Algo: G-APM} with $\gamma=\frac{2}{3}$, this paper & $1/L$ & $\mu/3L$ & Theorem~\ref{Thm: rate_SC_GAPM} \\ 
Monotone \eqref{Algo: M-APM} with $\gamma=1$ & $1/2L$ & $\mu/16L$ & \cite{Shi_2024_Monotone} \\
Monotone \eqref{Algo: M-APM} with $\gamma=1$, this paper & $1/2L$ & $\mu/4L$ & Corollary~\ref{Cor: rate_base_SC} \\
Monotone \eqref{Algo: M-APM} with $\gamma=\frac{2}{3}$, this paper & $1/L$ & $\mu/3L$ & Corollary~\ref{Cor: rate_fast_SC} \\
\bottomrule
\end{tabular}
\end{table}

The remainder of the paper is organized as follows. In Section \ref{Sec: energy_estimations}, we introduce an energy sequence and present some useful estimates. The convergence analysis for \eqref{Algo: M-APM} follows in Section \ref{Sec: convex} for convex functions, and in Section \ref{Sec: strongly convex} for strongly convex functions. In Section \ref{Sec: extensions}, we propose a variant of \eqref{Algo: APM} and establish a linear convergence rate under strong convexity. We then conclude in Section \ref{Sec: conclusions}.

\section{Energy estimations}\label{Sec: energy_estimations}
In this section, we introduce an energy sequence and derive its decrease along the sequences generated by \eqref{Algo: M-APM}. These results will be used in the forthcoming convergence analysis.

\subsection{Preliminaries}
To facilitate the convergence analysis, we define
\begin{equation}\label{E: Gs}
G_s(x) := \frac{ x - \prox_{sg}\left( x-s\nabla f(x) \right) }{s}.
\end{equation}
As shown in \cite[Lemma 7]{Bao_2023}, if $f$ is $\mu$-strongly convex, $s\in\left(0,\frac{1}{L}\right]$ and $x,y\in H$, we have 
\begin{equation} \label{E: prox_bound_SC}
F\left(x-sG_s(x)\right)  
  \le F(y) + \left\langle G_s(x), x-y \right\rangle
      - \frac{s(2 - sL)}{2} \| G_s(x) \|^2 - \frac{\mu}{2}\| x-y \|^2.    
\end{equation}
If $f$ is just convex, the inequality holds with $\mu=0$.

Using \eqref{E: Gs}, we have $z_k = x_k - sG_s(x_k)$, and \eqref{Algo: M-APM} gives:
\begin{equation}\label{E: iterates}
(t_{k+1}-1)(x_{k+1}-y_{k+1}) - (t_k-1)(x_k-y_k) + (x_{k+1}-x_k) 
= -\gamma s t_k G_s(x_k).
\end{equation}

Throughout the paper, the following standing assumption holds:
\begin{hypothesis}\label{Hypo: F}
Let $F=f+g$, where $f:H\to\R$ is $\mu$-strongly convex and $L$-smooth for $L\ge\mu\ge 0$, and $g:H\to\R\cup\{\infty\}$ is convex, proper and lower-semicontinuous.
\end{hypothesis}

\begin{remark}
If $f$ is only convex, Hypothesis \ref{Hypo: F} holds with $\mu=0$.
\end{remark}

\subsection{Energy decrease}
Our convergence analysis relies on the energy sequence $(E_k)_{k\ge 0}$, given by
\begin{equation}\label{E: E_k}
E_k=E_k(x^*) := \frac{1}{2}\| \xi_k - x^* \|^2 + \gamma s t_k (t_k-1) \left( F(y_k) - F^* \right),
\end{equation}
where $\xi_k = (t_k-1)(x_k-y_k) + x_k$.

\begin{remark}\label{Rem: E_k_compare}
In spite of the resemblence, the energy sequence $E_k$ is different from those of \cite{Beck_2009_Monotone,Shi_2024_Monotone} (with $\gamma=1$ in mind), whose $\xi_k$ takes the form
$$ \xi_k = t_k(z_k-y_k) + y_k = (t_k-1)(x_k-y_k) - st_kG_s(x_k) + x_k. $$
\end{remark} 

Using \eqref{E: E_k}, we have the following:

\begin{proposition}\label{Prop: E_k_diff_bound}
Let Hypothesis \ref{Hypo: F} hold. Let $(x_k)_{k\ge 1}$, $(y_k)_{k\ge 1}$ and $(z_k)_{k\ge 1}$ be generated by \eqref{Algo: M-APM} with $s\in \left( 0, 1/L \right]$ and $\gamma\in(0,\infty)$, and consider the sequence $(E_k)_{k\ge 0}$ defined by \eqref{E: E_k}. Then,
\begin{align*}
E_{k+1} - E_k
&\le - \frac{\gamma t_k^2(2-\gamma-sL)}{2}\| sG_s(x_k) \|^2
   - \frac{\mu \gamma s t_k(t_k-1)}{2} \|x_k-y_k \|^2 \\ 
&\quad - \frac{\mu\gamma s t_k}{2}\|x_k-x^* \|^2.
\end{align*}
\end{proposition}

\begin{proof}
By \eqref{E: E_k}, we have
\begin{equation}\label{E: E_k_diff}
\begin{aligned}
E_{k+1} - E_k
&= \left( \frac{1}{2}\| \xi_{k+1} - x^* \|^2 - \frac{1}{2}\| \xi_k - x^* \|^2 \right)
 + \gamma s t_{k+1} (t_{k+1}-1) \left( F(y_{k+1}) - F^* \right) \\
&\quad - \gamma s t_k(t_k-1) \left( F(y_k) - F^* \right).
\end{aligned}
\end{equation}
On the other hand, \eqref{E: iterates} implies that
$$\xi_{k+1} - \xi_k = -\gamma s t_k G_s(x_k),$$
so that
$$ \| \xi_{k+1} - \xi_k \|^2 = \gamma^2 t_k^2 \| s G_s(x_k) \|^2, $$
and
\begin{align*}
\langle \xi_k-x^*, \xi_{k+1} - \xi_k \rangle 
&= \langle (t_k-1)(x_k-y_k) + (x_k-x^*), -\gamma s t_k G_s(x_k) \rangle\\
&= - \gamma t_k (t_k-1) \langle sG_s(x_k), x_k-y_k \rangle 
 - \gamma t_k \langle s G_s(x_k), x_k-x^* \rangle.  
\end{align*}
Since
$$\frac{1}{2}\| \xi_{k+1} - x^* \|^2 - \frac{1}{2}\| \xi_k - x^* \|^2 = \langle \xi_k - x^*, \xi_{k+1} - \xi_k \rangle + \frac{1}{2}\| \xi_{k+1} - \xi_k \|^2,$$
we arrive at
\begin{equation}\label{E: phi_k_diff}
\begin{aligned}
\frac{1}{2}\| \xi_{k+1} - x^* \|^2 - \frac{1}{2}\| \xi_k - x^* \|^2 
&= - \gamma t_k (t_k-1) \langle sG_s(x_k), x_k-y_k \rangle \\ 
&\quad - \gamma t_k \langle s G_s(x_k), x_k-x^* \rangle
   + \frac{\gamma^2 t_k^2}{2}\| sG_s(x_k) \|^2.
\end{aligned}
\end{equation}
Noting that $F(y_{k+1}) \le F(z_k)$, and using \eqref{E: prox_bound_SC} with $x=x_k$ and $y=y_k$, we obtain
$$ F(y_{k+1}) - F(y_k) 
\le \langle G_s(x_k), x_k - y_k \rangle 
  - \frac{s(2-sL)}{2}\| G_s(x_k) \|^2  
  - \frac{\mu}{2}\|x_k-y_k \|^2. $$
Likewise, setting $x=x_k$ and $y=x^*$ in \eqref{E: prox_bound_SC}, we obtain
\begin{equation}\label{E:DeltaE_1}
F(y_{k+1}) - F^*
\le \langle G_s(x_k), x_k-x^* \rangle    
    - \frac{s(2-sL)}{2}\| G_s(x_k) \|^2
    - \frac{\mu}{2}\|x_k-x^* \|^2.
\end{equation}
Using these two inequalities, together with \eqref{E: phi_k_diff}, in \eqref{E: E_k_diff}, it follows that
\begin{align*}
E_{k+1} - E_k  
&\le - \gamma s\left[ t_k^2 - t_{k+1}(t_{k+1}-1) \right] \left( F(y_{k+1}) - F^* \right) 
  - \frac{\gamma t_k^2(2-\gamma-sL)}{2}\| sG_s(x_k) \|^2 \\  
&\quad - \frac{\mu \gamma s t_k(t_k-1)}{2}\|x_k-y_k \|^2
  - \frac{\mu\gamma s t_k}{2}\|x_k-x^* \|^2 \\
&= - \frac{\gamma t_k^2(2-\gamma-sL)}{2}\| sG_s(x_k) \|^2
   - \frac{\mu \gamma s t_k(t_k-1)}{2} \|x_k-y_k \|^2 \\
&\quad - \frac{\mu\gamma s t_k}{2}\|x_k-x^* \|^2,
\end{align*}
in view of $t_{k+1}(t_{k+1}-1)=t_k^2$.
\end{proof}

\begin{remark}\label{Rem: E_k_decreasing}
In particular, if $\gamma\in(0,2-sL]$, the sequence $(E_k)_{k\ge 0}$ is nonincreasing.
\end{remark}

\section{Convergence analysis I: the convex case}\label{Sec: convex}
\begin{theorem}\label{Thm: convergence_cvx}
Let Hypothesis \ref{Hypo: F} hold with $\mu=0$. Let $(x_k)_{k\ge 1}$, $(y_k)_{k\ge 1}$ and $(z_k)_{k\ge 1}$ be generated by \eqref{Algo: M-APM} with 
$$s\in \left( 0, 1/L \right] \qbox{and} \gamma\in(0,2-sL].  $$
Then for every $k\ge 1$, we have
$$ F(y_k) - F^* \le \frac{\| x_0 - x^* \|^2}{2\gamma st_k(t_k-1)}. $$
Besides, every weak subsequential limit point of $y_k$, as $k\to\infty$, minimizes $F$.
\end{theorem}

\begin{proof}
Setting $\mu=0$ in Proposition \ref{Prop: E_k_diff_bound}, we obtain $E_{k+1}-E_k\le 0$, since $2-\gamma-sL\ge 0$. It follows that
$$ F(y_k) - F^* \le \frac{E_k}{\gamma st_k(t_k-1)} \le \frac{E_0}{\gamma st_k(t_k-1)},\ \forall\,k\ge 1. $$
With $t_0=1$ in mind, we obtain from \eqref{E: E_k} that
$$ E_0 = \frac{1}{2}\| x_0 - x^* \|^2. $$
This gives
$$ F(y_{k}) - F^* \le \frac{\| x_0 - x^* \|^2}{2\gamma st_k(t_k-1)},\ \forall\,k\ge 1. $$
Observing that $F$ is weakly lower-semicontinuous, and that $t_k\to\infty$ as $k\to\infty$, we obtain the minimizing property.
\end{proof}

\begin{remark}
A similar but slightly better convergence rate $\mathcal{O}\left( \frac{1}{k^2} \right)$ has been reported in \cite{Beck_2009_Monotone} for the case $\gamma=1$. The improvement is due to the different energy sequences used for the proof, as discussed in Remark \ref{Rem: E_k_compare}. However, the energy sequence introduced in \cite{Beck_2009_Monotone} seems not applicable to the convergence proof for \eqref{Algo: M-APM} in the presence of the correction term \eqref{E: correction_term}.
\end{remark}

In what follows, we establish the boundedness of the iterates generated by \eqref{Algo: M-APM}.

\begin{lemma}\label{Lem: bounded_iterates}
Let Hypothesis \ref{Hypo: F} hold with $\mu=0$. Let $(x_k)_{k\ge 1}$, $(y_k)_{k\ge 1}$ and $(z_k)_{k\ge 1}$ be generated by \eqref{Algo: M-APM} with $s\in \left( 0, 1/L \right]$ and $\gamma\in(0,2-sL]$. Then, $(x_k)_{k\ge 0}$, $(y_k)_{k\ge 0}$ and $(z_k)_{k\ge 0}$ are bounded. 
\end{lemma}

\begin{proof}
Consider the sequence $(E_k)_{k\ge 0}$ defined by \eqref{E: E_k}, where $x^*$ is an arbitrary minimizer of $F$. By Proposition \ref{Prop: E_k_diff_bound}, $E_{k+1}\le E_k$. It follows that
$$ \| \xi_k - x^* \|^2 \le 2E_k \le 2E_0, $$
which means $(\xi_k)_{k\ge 0}$ is bounded. Let $\| \xi_k \| \le M$ for some constant $M>0$. 
Recalling that
$$ \xi_k = (t_k-1)(x_k-y_k) + x_k, $$
we have
$$ x_k = \left( 1 - \frac{1}{t_k} \right)y_k + \frac{1}{t_k}\xi_k, $$
so that
\begin{equation}\label{E: x_k_bound}
\| x_k \| \le \left( 1 - \frac{1}{t_k} \right)\| y_k \| + \frac{1}{t_k}M.
\end{equation}
Using \eqref{E: iterates}, we have
\begin{equation}\label{E: xi_k_diff}
\xi_{k+1} - \xi_k = \gamma t_k(z_k-x_k),
\end{equation}
so that
$$ z_k = x_k + \frac{1}{\gamma t_k}(\xi_{k+1}-\xi_k). $$
This gives
$$ \| z_k \| 
\le \| x_k \| + \frac{2M}{\gamma t_k} 
\stackrel{\eqref{E: x_k_bound}}{\le} 
\left( 1 - \frac{1}{t_k} \right)\| y_k \| + \frac{1}{t_k}\left( 1 + \frac{2}{\gamma} \right)M
\le \max\left\{ \| y_k \|, \left( 1 + \frac{2}{\gamma} \right)M \right\}. $$
From the iterate for $y_{k+1}$ in \eqref{Algo: M-APM}, we have
$$\| y_{k+1} \| 
\le \max\left\{ \| z_k \|, \| y_k \| \right\} 
\le \max\left\{ \| y_k \|, \left( 1 + \frac{2}{\gamma} \right)M \right\}
\le \max\left\{ \| y_0 \|, \left( 1 + \frac{2}{\gamma} \right)M \right\},$$
so that $(y_k)_{k\ge 0}$ is bounded. Then, $(z_k)_{k\ge 0}$ is bounded too. In view of \eqref{E: x_k_bound}, we obtain the boundedness of $(x_k)_{k\ge 0}$.
\end{proof}

The following result \cite[Lemma A.4]{Radu_2025} will be useful in the sequel:

\begin{lemma}\label{Lem: convergence_limit}
Let $(u_k)_{k\ge 0}$ be a real sequence and $(\zeta_k)_{k\ge 0}$ be positive such that $\sum_{k=0}^{\infty}\frac{1}{\zeta_k}=\infty$. If
$ \lim_{k\to\infty}[ u_{k+1} + \zeta_k( u_{k+1} - u_k ) ] = b\in\R, $
then,
$ \lim_{k\to\infty}u_k = b. $ 
\end{lemma}

Now we are ready to prove that all weak subsequential limit points of the iterates generated by \eqref{Algo: M-APM} are minimizers of $F$.

\begin{theorem}
Let Hypothesis \ref{Hypo: F} hold with $\mu=0$. Let $(x_k)_{k\ge 1}$, $(y_k)_{k\ge 1}$ and $(z_k)_{k\ge 1}$ be generated by \eqref{Algo: M-APM} with $s\in \left( 0, 1/L \right]$ and $\gamma\in(0,2-sL)$. Then, all weak subsequential limit points of $(x_k)_{k\ge 0}$, $(y_k)_{k\ge 0}$ and $(z_k)_{k\ge 0}$ belong to $\argmin(F)$. 
\end{theorem}

\begin{proof}
From Lemma \ref{Lem: bounded_iterates}, the sequences $(x_k)_{k\ge 0}$, $(y_k)_{k\ge 0}$ and $(z_k)_{k\ge 0}$ are bounded, so that they each have at least one weak subsequential limit point. From Theorem \ref{Thm: convergence_cvx}, all weak subsequential limit points of $(y_k)_{k\ge 0}$ are minimizers of $F$. Setting $\mu=0$ in Proposition \ref{Prop: E_k_diff_bound} and with $\gamma\in(0,2-sL)$ in mind, we have
$$ E_{k+1} - E_k \le - \frac{\gamma t_k^2(2-\gamma-sL)}{2}\| sG_s(x_k) \|^2, $$
so that
$$ \frac{\gamma(2-\gamma-sL)}{2}\sum_{k=0}^\infty t_k^2 \| sG_s(x_k) \|^2 \le E_0 < \infty. $$
This gives
\begin{equation}\label{E: G_s_limit}
\lim_{k\to\infty}\| s t_k G_s(x_k) \| = 0.
\end{equation} 
Since $sG_s(x_k) = x_k-z_k$, we have
$$ \lim_{k\to\infty} \| t_k(x_k-z_k) \| = 0, $$
which implies that $(x_k)_{k\ge 0}$ shares the same weak subsequential limit points with $(z_k)_{k\ge 0}$. Using \eqref{E: G_s_limit} in \eqref{E: iterates}, we obtain
$$ \lim_{k\to\infty} \left\| t_{k+1}\left[ x_{k+1}-\left( 1 - \frac{1}{t_{k+1}} \right)y_{k+1} \right] - t_k\left[ x_k-\left( 1 - \frac{1}{t_k} \right)y_k \right] \right\| = 0, $$
which implies that
$$ \lim_{k\to\infty} (t_{k+1}-t_k)\| u_{k+1} \| + t_k\left( \| u_{k+1} \| - \| u_k \| \right) = 0, $$
where
$$ u_k := x_k-\left( 1 - \frac{1}{t_k} \right)y_k.$$
It follows that
$$ \lim_{k\to\infty} \| u_{k+1} \| + \frac{t_k}{t_{k+1}-t_k}\left( \| u_{k+1} \| - \| u_k \| \right) = 0. $$
By invoking Lemma \ref{Lem: convergence_limit}, we deduce that $\lim_{k\to\infty} \| u_k \| = 0$, and then
$$ \lim_{k\to\infty}\| x_k-y_k \| = 0. $$
This implies that $(x_k)_{k\ge 0}$ shares the same weak subsequential limit points with $(y_{k})_{k\ge 0}$, and allows us to conclude the proof.
\end{proof}

\begin{remark}
Due to existence of the $(z_k)$ sequence in \eqref{Algo: M-APM}, it seems difficult to gain a weak convergence of the iterates, which is in contrast to the analysis in \cite{Ryu_2025,Radu_2025_iterate} for the FISTA algorithm \eqref{Algo: APM}.
\end{remark}

\section{Convergence analysis II: the strongly convex case}\label{Sec: strongly convex}

\begin{proposition} \label{Prop: Bound_Ek+1}
Let Hypothesis \ref{Hypo: F} hold with $\mu>0$. Consider \eqref{Algo: M-APM} with $s\in \left( 0, 1/L \right]$ and $\gamma\in(0,2-sL)$. Consider the sequence $(E_k)_{k\ge 0}$ defined by \eqref{E: E_k}. Then,
$$E_{k+1}
\le (t_k-1)^2 \| x_k-y_k \|^2 
 + 2 \| x_k-x^* \|^2 
 + \left[ 1- \mu s(2-4\gamma-sL) \right] \frac{\gamma t_k^2}{2\mu s} \| s G_s(x_k) \|^2.$$
\end{proposition}

\begin{proof}
By definition of $\xi_k$, we have
$$ \xi_{k+1} 
= (t_{k+1}-1)(x_{k+1}-y_{k+1}) + x_{k+1}
\stackrel{\eqref{E: iterates}}{=} (t_k-1)(x_k-y_k) + x_k - \gamma st_k G_s(x_k). $$
Therefore,
\begin{equation}\label{E: phi_k_bound_SC}
\frac{1}{2}\| \xi_{k+1} - x^* \|^2
\le (t_k-1)^2 \| x_k-y_k \|^2 
 + 2 \| x_k-x^* \|^2 \\
 + 2 \gamma^2 t_k^2 \| sG_s(x_k) \|^2.
\end{equation}
On the other hand, \eqref{E:DeltaE_1} gives
\begin{align*}
F(y_{k+1}) - F^* 
&\le \langle G_s(x_k), x_k-x^* \rangle - \frac{\mu}{2}\|x_k-x^* \|^2
    - \frac{s(2-sL)}{2}\| G_s(x_k) \|^2 \\
&\le \left[ \frac{1- \mu s(2-sL) }{2\mu} \right] \| G_s(x_k) \|^2.
\end{align*}
With this and \eqref{E: phi_k_bound_SC} in mind, and recalling that 
\begin{align*}
E_{k+1} 
&= \frac{1}{2}\| \xi_{k+1} - x^* \|^2 + \gamma s t_{k+1}(t_{k+1}-1)\left( F(y_{k+1}) - F^*\right) \\
&= \frac{1}{2}\| \xi_{k+1} - x^* \|^2 + \gamma s t_k^2 \left( F(y_{k+1}) - F^*\right),
\end{align*}
we obtain the desired result.
\end{proof}

\begin{remark}
This bound is obtained following the analysis in \cite{Bao_2023}.
\end{remark}

Now we are ready to derive the linear convergence rate of the function values for \eqref{Algo: M-APM} under strong convexity.

\begin{theorem}\label{Thm: rate_SC}
Let Hypothesis \ref{Hypo: F} hold with $\mu>0$. Let $(x_k)_{k\ge 1}$, $(y_k)_{k\ge 1}$ and $(z_k)_{k\ge 1}$ be generated by \eqref{Algo: M-APM} with 
$$s\in \left( 0, \frac{1}{L} \right] \qbox{and} \gamma\in\left( \frac{1-sL}{4},2-sL \right).  $$
Then for every $k\ge 2$, we have 
$$ F(y_k) - F^* \le \frac{\| x_0 - x^* \|^2}{2\gamma st_k(t_k-1)}(1+\rho)^{-k+2},$$
where
$$ \rho = \min\left\{ \frac{\mu s(2-\gamma-sL)}{ 1- \mu s(2-4\gamma-sL) },\ 
\frac{\mu\gamma s}{2}
 \right\}. $$
\end{theorem}

\begin{proof}
By comparing Propositions \ref{Prop: E_k_diff_bound} and \ref{Prop: Bound_Ek+1}, we deduce that
$$ E_{k+1} - E_k \le -r_k E_{k+1}, $$
where
\begin{equation}\label{E: r_k}
r_k = \min\left\{ \frac{\mu s(2-\gamma-sL)}{ 1- \mu s(2-4\gamma-sL) },\ 
\frac{\mu\gamma st_k}{2(t_k-1)},\ 
\frac{\mu\gamma s t_k}{4}
 \right\}.
\end{equation}
Recalling the definition of $(t_k)$ in \eqref{E: t_k}, we have $t_k\ge 2$ for every $k\ge 2$. Whence, we deduce that
$$ E_{k+1} - E_k \le -\rho E_{k+1},\quad \forall k\ge 2, $$
where
$$ \rho = \min\left\{ \frac{\mu s(2-\gamma-sL)}{ 1- \mu s(2-4\gamma-sL) },\ 
\frac{\mu\gamma s}{2}
 \right\}. $$
It follows that
$$ E_k \le E_2(1+\rho)^{-(k-2)} \le E_0(1+\rho)^{-(k-2)},\quad \forall k\ge 2, $$
in view of Remark \ref{Rem: E_k_decreasing}. Then, we arrive at 
$$ F(y_k) - F^* \le \frac{E_k}{\gamma st_k(t_k-1)} \le \frac{E_0}{\gamma st_k(t_k-1)}(1+\rho)^{-(k-2)},\quad\forall k\ge 2. $$
The result follows by recalling that $E_0 = \frac{1}{2}\| x_0 - x^* \|^2$.
\end{proof}

\begin{remark}
We use the convention of $\frac{c}{0}=\infty$ for any $c>0$ so that the second term in the right-hand side of \eqref{E: r_k} is well defined, in view of $t_0=1$. 
\end{remark}

\begin{remark}\label{Rem: rate_base_SC}
In case $\gamma=1$, the step size should satisfy $s\in(0,1/L)$ for a linear convergence guarantee. The critical step size $s=1/L$ is not allowed. Note that
$$ \rho = \min\left\{ \frac{\mu s(1-sL)}{ 1+ \mu s(2+sL) },\ 
\frac{\mu s}{2}
 \right\} 
 \sim 
 \min\left\{ \mu s(1-sL),\ 
 \frac{\mu s}{2} \right\}, $$
assuming $\mu s \ll 1$. As one can see, $\rho$ is maximized with $s=1/(2L)$ and $\rho_{\max} = \frac{\mu}{4L+5\mu}$.
\end{remark}

\begin{remark}\label{Rem: rate_fast_SC}
In case $s=1/L$, the value of $\gamma$ should satisfy $\gamma\in(0,1)$ so as to obtain a linear convergence. Note that
$$ \rho = \min\left\{ \frac{\mu(1-\gamma)}{ L- \mu (1-4\gamma) },\ 
\frac{\mu\gamma}{2L}
 \right\}
 \sim
 \frac{\mu}{L}\min\left\{ (1-\gamma), \gamma/2 \right\}, $$
assuming $\mu/L\ll 1$. It is easy to show that $\rho$ is maximized with $\gamma=\frac{2}{3}$ and $\rho_{\max} = \frac{\mu}{3L+5\mu}$.
\end{remark}

Based on Remarks \ref{Rem: rate_base_SC} and \ref{Rem: rate_fast_SC}, we obtain the following results.

\begin{corollary}\label{Cor: rate_base_SC}
Let Hypothesis \ref{Hypo: F} hold with $\mu>0$. Let $(x_k)_{k\ge 1}$, $(y_k)_{k\ge 1}$ and $(z_k)_{k\ge 1}$ be generated by \eqref{Algo: M-APM} with 
$$s=\frac{1}{2L} \qbox{and} \gamma=1.  $$
Then for every $k\ge 2$, we have 
$$ F(y_k) - F^* \le \frac{L\| x_0 - x^* \|^2}{t_k(t_k-1)}(1+\rho)^{-k+2},\qbox{with}
\rho = \frac{\mu}{4L+5\mu}.$$
\end{corollary} 

\begin{remark}
This rate is faster than that of \cite{Shi_2024_Monotone}, where the best rate is $\rho\sim \frac{\mu}{16L}$ for $s=\frac{1}{2L}$. The presented proof follows the idea of \cite{Bao_2023}, but one should also notice that the energy sequence $E_k$ in \eqref{E: E_k} is different from the one in \cite{Bao_2023}, as discussed in Remark \ref{Rem: E_k_compare}. We introduce a different $E_k$ to accommodate the correction term \eqref{E: correction_term} in the proof, which will contribute to a faster convergence rate in what follows.  
\end{remark}

\begin{corollary}\label{Cor: rate_fast_SC}
Let Hypothesis \ref{Hypo: F} hold with $\mu>0$. Let $(x_k)_{k\ge 1}$, $(y_k)_{k\ge 1}$ and $(z_k)_{k\ge 1}$ be generated by \eqref{Algo: M-APM} with 
$$s=\frac{1}{L} \qbox{and} \gamma=\frac{2}{3}.  $$
Then for every $k\ge 2$, we have 
$$ F(y_k) - F^* \le \frac{3L\| x_0 - x^* \|^2}{4 t_k(t_k-1)}(1+\rho)^{-k+2},\qbox{with}\rho = \frac{\mu}{3L+5\mu}.$$
\end{corollary}

\begin{remark}
This rate is even faster than that of Corollary \ref{Cor: rate_base_SC}, which is attributed to the introduction of $\gamma$ in \eqref{Algo: M-APM}. Note also that the critical step size $s=1/L$ is allowed, in contrast to the case with $\gamma=1$.  
\end{remark}

\section{Extensions to accelerated proximal methods}\label{Sec: extensions}
In view of the faster convergence rate presented in Corollary \ref{Cor: rate_fast_SC} for \eqref{Algo: M-APM}, we propose the following generalized accelerated proximal method:
\begin{equation}\label{Algo: G-APM}\tag{G-APM}
\left\{
\begin{aligned}
y_{k+1} &= \prox_{sg}\left( x_k - s\nabla f(x_k) \right),\\
x_{k+1} &= y_{k+1} + \frac{t_k-1}{t_{k+1}}(y_{k+1}-y_k) + \frac{t_k}{t_{k+1}}(\gamma-1)(y_{k+1}-x_k),
\end{aligned}
\right.
\end{equation}
where $s\in \left( 0, 1/L \right]$, $\gamma\in(0,\infty)$ and $(t_k)$ is given by \eqref{E: t_k}. When $\gamma=1$, \eqref{Algo: G-APM} reduces to \eqref{Algo: APM}, which enjoys a linear convergence rate $\mathcal{O}\left( \frac{1}{k^2}(1+\rho)^{-k}\right)$, with $\rho=\frac{\mu}{4L}$ for $s=\frac{1}{2L}$ in case $f$ is $\mu$-strongly convex \cite{Bao_2023}. In this section, we will show that the constant $\rho$ can be further improved to $\rho=\frac{\mu}{3L}$, when $\gamma=\frac{2}{3}$ and $s=\frac{1}{L}$.

Using \eqref{E: Gs}, we have $y_{k+1} = x_k - sG_s(x_k)$. This, together with \eqref{Algo: G-APM}, gives the iterates:
\begin{equation}\label{E: iterates_APM}
(t_{k+1}-1)(x_{k+1}-y_{k+1}) - (t_k-1)(x_k-y_k) + (x_{k+1}-x_k) = -\gamma s t_k G_s(x_k).
\end{equation} 
Observe that \eqref{E: iterates_APM} is identical to \eqref{E: iterates} in its form. This implies that Theorem \ref{Thm: rate_SC} and Corollary \ref{Cor: rate_fast_SC} also hold for \eqref{Algo: G-APM}. To avoid repetition, we give the following result of interest and leave the proof to the reader: 

\begin{theorem}\label{Thm: rate_SC_GAPM}
Let Hypothesis \ref{Hypo: F} hold with $\mu>0$. Let $(x_k)_{k\ge 1}$, $(y_k)_{k\ge 1}$ and $(z_k)_{k\ge 1}$ be generated by \eqref{Algo: G-APM} with 
$$s=\frac{1}{L} \qbox{and} \gamma=\frac{2}{3}.  $$
Then for every $k\ge 2$, we have 
$$ F(y_k) - F^* \le \frac{3L\| x_0 - x^* \|^2}{4 t_k(t_k-1)}(1+\rho)^{-k+2},\qbox{with}\rho = \frac{\mu}{3L+5\mu}.$$
\end{theorem}

\begin{remark}
This linear rate is faster than those of \cite{Shi_2024,Bao_2023} for \eqref{Algo: APM}, where the best rate is $\mathcal{O}\left( \frac{1}{k^2}(1+\rho)^{-k}\right)$ with $\rho=\frac{\mu}{4L}$ when $s=\frac{1}{2L}$. Note also that the critical step size $s=\frac{1}{L}$ is allowed in our result, thanks to introduction of the term involving $\gamma$. For \eqref{Algo: APM} with $s=\frac{1}{L}$, a linear rate $\mathcal{O}\left( \left(1 - \frac{\mu^2}{4L^2} \right)^{k} \right)$ was derived in \cite{Bao_2023}, which is much slower than ours.   
\end{remark}

\section{Conclusions}\label{Sec: conclusions}
We analyze the convergence rate for \eqref{Algo: M-APM}, and establish a fast rate $\mathcal{O}\left( \frac{1}{k^2} \right)$ for the function values in the convex setting, and a linear rate $\mathcal{O}\left( \frac{1}{k^2}\left( 1 + \frac{\mu}{3L} \right)^{-k} \right)$ under strong convexity. In the convex setting, we also prove boundedness of the iterates, and show that their weak subsequential limit points are all minimizers of the objective function. We propose a generalized variant of \eqref{Algo: APM} and obtain a faster linear convergence rate under strong convexity.

\begin{acknowledgements}
This work was partially funded by the China Scholarship Council~202208520010, and also benefited from the support of the FMJH Program Gaspard Monge for optimization and operations research and their interactions with data science.
\end{acknowledgements}

\bibliography{myrefs}

\end{document}